\documentclass[12pt]{amsart}
\usepackage{amsmath,amsthm}

\def\pmatrix{\left(\begin{matrix}}
\def\endpmatrix{\end{matrix}\right)}

\def\Z{{\mathbb Z}}
\def\C{{\mathbb C}}

\def\cal{\mathcal}

\def\dd{{\rm d}}

\def\de{\delta}

\def\p{\partial}

\def\t{\theta}
\def\s{\sigma}
\def\T{\Theta}
\def\e{\varepsilon}

\def\a{\alpha}
\def\b{\beta}

\def\H{{\mathcal H}}

\def\tt#1#2{{\t\left[\begin{matrix}{#1}\\ {#2}\end{matrix}\right]}}

\theoremstyle{plain}
\newtheorem{thm}{Theorem}%[section]
\newtheorem{lm}[thm]{Lemma}
\newtheorem{prop}[thm]{Proposition}

\theoremstyle{definition}

\newtheorem{rem}[thm]{Remark}

\title{Two generalizations of Jacobi's derivative formula}
\author{Samuel Grushevsky, Riccardo Salvati Manni}
\address{Mathematics Department, Princeton University, Fine Hall,
Washington Road, Princeton, NJ 08544, USA}
\email{sam@math.princeton.edu}
\address{Dipartimento di Matematica, Universit\`a ``La Sapienza'',
Piazzale A. Moro 2, Roma, I 00185, Italy}
\email{salvati@mat.uniroma1.it}
\thanks{First author partially supported by NSF Mathematical Sciences
Postdoctoral Research Fellowship}
\begin{document}
\begin{abstract}
In this paper we generalize Jacobi's derivative formula,
considered as an identity for theta functions with characteristics
and their derivatives, to higher genus/dimension. By applying the
methods developed in our previous paper \cite{gsm}, several
generalizations to Siegel modular forms are obtained. These
generalizations are identities satisfied by theta functions with
characteristics and their derivatives at zero. Equating all the
coefficients of the Fourier expansion of these relations to zero
yields non-trivial combinatorial identities.
\end{abstract}

\maketitle

\section{Definitions and notations}
We denote by $\H_g$ the {\it Siegel upper half-space} --- the
space of symmetric complex $g\times g$ matrices with positive
definite imaginary part. For $\e,\de\in (\Z/2\Z)^g$, thought of as
vectors of zeros and ones, $\tau\in\H_g$ and $z\in \C^g$, the {\it
theta function with characteristic $[\e,\de]$} is
$$
\tt\e\de(\tau,z):=\sum\limits_{m\in\Z^g} \exp \pi i \left[
^t(m+\frac{\e}{2})\tau(m+\frac{\e}{2})+2\ ^t(m+\frac{\e}{2})(z+
\frac{\de}{2})\right].
$$
A {\it characteristic} $[\e,\de]$ is called {\it even} or {\it
odd} depending on whether it is even or odd as a function of $z$,
which corresponds to the scalar product $\e\cdot\de\in\Z/2\Z$
being zero or one, respectively. The number of even (resp. odd)
theta characteristics is $2^{g-1}(2^g+1)$ (resp.
$2^{g-1}(2^g-1)$). For $\e\in(\Z/2\Z)^g$ the {\it second order
theta function} with characteristic $\e$ is
$$
\T[\e](\tau,z):=\tt{\e}{0}(2\tau,2z).
$$

A {\it theta constant} is the evaluation at $z=0$ of a theta
function. We drop the argument $z=0$ in the notations for theta
constants. Obviously all odd theta constants vanish identically,
and thus there are $2^{g-1}(2^g+1)$ non-trivial theta constants
with characteristics, and $2^g$ theta constants of the second
order.

A triplet of characteristics $[\e_1,\de_1]$, $[\e_2,\de_2]$,
$[\e_3,\de_3]$ is called {\it azygetic} if
$$
(-1)^{\e_1\cdot\de_1+\e_2\cdot\de_2+\e_3\cdot\de_3+
(\e_1+\e_2+\e_3)\cdot(\de_1+\de_2+\de_3)}=-1.
$$
A sequence of $2g+2$ characteristics $[\e_1,\de_1],\dots,
[\e_{2g+2},\de_{2g+2}]$ is called a {\it special fundamental
system} if the first $g$ characteristics are odd, the remaining
are even and any triple of characteristics in it is azygetic.

In the genus $1$ case one of the main identities for theta
functions is Jacobi's derivative formula:
\begin{equation}
\label{jac0}
\frac{\dd}{\dd z}\tt11(\tau,z)|_{z=0}=-\pi\tt00(\tau)
\tt10(\tau) \tt01(\tau).
\end{equation}
There is a long history of possible generalizations of this
formula to higher genus. We consider $g$ odd characteristics
$[\e_1,\de_1],\dots, [\e_g,\de_g]$, and define their {\it jacobian
determinant} to be
\begin{equation}
\label{jac}
\begin{matrix}
D([\e_1,\de_1],\dots[\e_g,\de_g])(\tau):=\qquad\qquad\hfill\\
\qquad\hfill \pi^{-g} grad\,\,\tt{\e_{1}}{\de_{1}} \wedge
grad\,\,\tt{\e_{2}}{\de_{2 }}\wedge \dots \wedge
grad\,\,\tt{\e_{g}}{\de_{g}}(\tau, 0).
\end{matrix}
\end{equation}
Essentially the problem of generalizing Jacobi's derivative
formula consists in expressing some linear combinations of
jacobian determinants of $g$ distinct odd theta functions as
polynomials or rational functions in theta constants. When $g=2$,
such formulas are due to Rosenhain, Frobenius, Thomae, Fay, Igusa:
we refer to \cite{thomae}, \cite{fr885}, \cite{fay}, \cite{ig2}
and \cite{sm3} for exact statements and a precise description of
the situation. We recall from these works that there is a precise
conjectural formula, which has been proven for $g\leq 5$.
Moreover, for $g\leq 3$ the equality
\begin{equation}
\label{jac1} D([\e_1,\de_1],\dots,[\e_g,\de_g])(\tau)=\pm
\tt{\e_{g+1}}{\de_{g+1}}(\tau)\dots \tt{\e_{2g+2}}{\de_{2g+2}}(\tau)
\end{equation}
holds if and only if the $2g+2$ characteristics appearing in it form
a special fundamental system.

Differential equations for genus $2$ theta constants have also
been studied by Ohyama \cite{Oh} and Zudilin \cite{Zu}; Grant
\cite{gr} obtains a nice relation involving only one partial
derivative.

Generalizations of Jacobi's derivative formula in another
direction, to higher level theta constants in one variable, are
derived and discussed in \cite{farkaskra} --- generalizing these
to the higher genus would also be very interesting.

A different generalization of Jacobi's derivative formula involves
higher order derivatives of theta functions. For example it makes
sense in genus $1$ to ask for the expression of
\begin{equation}
\label{jac2}
\det\pmatrix\tt00(\tau)\,\quad&\tt10(\tau )\,\quad\\
\frac{\dd^2}{\dd^2z}\tt00
(\tau,z)|_{z=0}&\frac{\dd^2}{\dd^2z}\tt10(\tau,z)|_{z=0}
\endpmatrix
\end{equation}
$$
=4\pi i\det\pmatrix\quad\tt00(\tau)&\quad\tt10(\tau)\\
\frac{\dd}{\dd\tau }\tt00(\tau)&\frac{\dd}{\dd\tau}\tt10(\tau)
\endpmatrix
$$
as a polynomial in theta constants and first-order derivatives
with respect to $z$ of theta functions, evaluated at $z=0$ (the
two determinants are equal by the heat equation). We know in fact
that
\begin{equation}
\label{jac3}
-\T[1](2\tau)^2\frac{\dd}{\dd\tau}(\T[0](\tau))/\T[1](\tau))=
\frac{i}{4\pi}\left(\frac{\dd}{\dd z}\tt11(\tau,z)|_{z=0}\right)^2
\end{equation}
or, equivalently,
\begin{equation}
\label{jac4}
\det\pmatrix\quad \T[0](\tau)&\quad \T[1](\tau )\\
\frac{\dd}{\dd\tau}
\T[0](\tau)&\frac{\dd}{\dd\tau}\T[1](\tau)\endpmatrix =
\frac{i}{4\pi}\left(\frac{\dd}{\dd
z}\tt11(\tau,z)|_{z=0}\right)^2.
\end{equation}
To prove this, one can invoke a modular argument, saying that both
sides are modular of the same weight, thus proportional, and the
constant can be easily computed. Of course such a proof is not very
revealing, and thus obtaining another proof would be desirable. As
explained to us by H. Farkas, this identity can also be deduced from
theorem 5.3 in chapter 2 of \cite{farkaskra} by applying Jacobi's
triple product identity, changing to the argument $\tau/2$ and then
verifying the resulting identity combinatorially for each
coefficient of the Fourier series.

In this paper we shall generalize both the identities (\ref{jac3})
and (\ref{jac4}) to higher genus. It would be interesting to
understand the combinatorial meaning of these generalizations
similarly to the one-variable identities above or obtain
alternative combinatorial proofs, but these questions lie beyond
the scope of the current work.

The main tool will be a consequence of Riemann's addition theorem
relating the first $z$-derivatives of odd theta functions with
characteristics to the second $z$-derivatives of second order
theta functions. This has also been the main tool in our paper
\cite{gsm}, where we showed that generically a principally
polarized abelian variety is uniquely determined by the gradients
of odd theta functions at $z=0$.

\begin{rem} We note that the classical generalization of Jacobi's
derivative formula can be given an interpretation in terms of
theta series with harmonic polynomial coefficients. In fact
monomials of degree $g+2$ in the theta constants are theta series
relative to the quadratic form $4\cdot 1_{g+2}$ and harmonic
polynomial ``1''. The jacobian determinants, on the other hand,
are theta series relative to the quadratic form $4\cdot 1_{g}$ and
harmonic polynomial ``$\det$'', \cite{Igu}.

This is the simplest pair of theta series with harmonic polynomial
coefficients. Our generalizations can also be interpreted  in this
way. For example in genus one, while, as we wrote, in the first
two cases the  harmonic polynomials are $1$ and $x$, in our case
it is a polynomial in two variables:  $x^2-y^2$.
\end{rem}

As a further consequence of our formulas, we shall give a
characterization of the locus of reducible principally polarized
abelian variety in terms of vanishing of certain derivatives of
odd theta functions.

\section{The symplectic group action}
Let $\Gamma_g:= Sp(2g,\Z)$ be the integral symplectic group; it
acts on $\H_g$ by
$$
M\cdot \tau:=(A\tau+B)(C \tau+D)^{-1},
$$
where $M=\pmatrix A&B\\ C&D\endpmatrix\in\Gamma_g.$ A period
matrix $\tau$ is called {\it reducible} if there exists
$M\in\Gamma_g$ such that
$$
M\cdot\tau=\pmatrix \tau_1&0\\
0&\tau_2\endpmatrix,\quad\tau_i\in\H_{g_i},\ g_1+g_2=g;
$$
otherwise we say that $\tau$ is irreducible.

We define the {\it level} subgroups of the symplectic group to be
$$
\Gamma_g(n):=\left\lbrace M=\pmatrix A&B\\ C&D\endpmatrix
\in\Gamma_g\, |\, M\equiv\pmatrix 1&0\\
0&1\endpmatrix\ {\rm mod}\ n\right\rbrace
$$
$$
\Gamma_g(n,2n):=\left\lbrace M\in\Gamma_g(n)\, |\, {\rm
diag}(A^tB)\equiv{\rm diag} (C^tD)\equiv0\ {\rm mod}\
2n\right\rbrace.
$$

A function $F:\H_g\to\C$ is called a {\it modular form of weight
$k$ with respect to $\Gamma\subset\Gamma_g$} if
$$
F(M\cdot\tau)=\det(C\tau+D)^kF(\tau),\quad \forall M \in\Gamma,\
\forall \tau\in\H_g
$$

The theta functions transform under the action of $\Gamma_g$ as
follows:
$$
\theta \bmatrix M\pmatrix \e\\ \de\endpmatrix\endbmatrix
(M\cdot\tau,\,^{t}(C\tau+D)^{-1}z)\qquad\qquad\qquad$$
$$\qquad\qquad\qquad=\phi(\e,\,\de,\,M,\,
\tau,\,z)\det(C\tau+D)^{\frac{1}{2}}\theta\bmatrix \e\cr \de
\endbmatrix(\tau,\,z),
$$
where
$$
M\pmatrix \e\cr \de\endpmatrix :=\pmatrix D&-C\cr
-B&A\endpmatrix\pmatrix \e\cr \de\endpmatrix+ \pmatrix {\rm
diag}(C \,^t D)\cr {\rm diag}(A\,^t B)\endpmatrix
$$
taken modulo 2, and $\phi(\e,\,\de,\,M,\,\tau,\,z)$ is some
complicated explicit function. For more details, we refer to
\cite{ig1} and \cite{farkasrauch}.

Theta constants with characteristics are modular forms of weight
$1/2$ with respect to $\Gamma_g(4,8)$. In this case $\phi(\e,\de,M
):=\phi(\e,\de,M,\tau,0)$ is an eighth root of unity that does not
depend on $\tau$.

Differentiating the theta transformation law above with respect to
some $z_i$ and then evaluating at $z=0$, we see that
$$
\frac{\p}{\p z_i}\theta\left[M\pmatrix\e\\ \de
\endpmatrix\right](M\cdot\tau,z)|_{z=0}\qquad\qquad\hfill
$$
$$
\hfill\quad=\phi(\e,\de,M)\det(C\tau+D)^{1/2}\sum\limits_j
(C\tau+D)_{ij}\frac{\p}{\p z_j}\tt\e\de( \tau,z)|_{z=0}.
$$
Denoting by $grad\,\,\tt\e\de( \tau)$ the gradient of the theta
function with respect to $z_1,z_2,\dots,z_g$ at $z=0$, the above
formula becomes
$$
grad\,\,\theta\!\left[M\pmatrix\e\\ \de
\endpmatrix\right]\!(M\cdot
\tau)\!=\!\phi(\e,\de,M)\det(C\tau+D)^{\frac{1}{2}}(C\tau+D)grad\,
\,\tt\e\de(\tau).
$$
As a consequence, the jacobian determinant
$D([\e_1,\,\de_1],\dots[\e_g,\,\de_g])(\tau)$ is a modular form of
weight $\frac{1}{2}g+1$ with respect to $\Gamma_g(4,8)$ (see
\cite{ig2} and \cite{sm3}).

\section{Some multilinear algebra}
For our purposes we need some results from linear algebra, which we
recall and prove for the sake of completeness. We are grateful to C.
De Concini, A. Maffei,  D. Zagier and one of the referees for useful
suggestions about these topics.

To any $A\in {\rm Mat}_{g\times g}(\C)$ we associate the
$(g-1)\times (g-1)$ matrix ${\tilde A}$ whose entries are the
determinants of $2\times 2$ minors of $A$ obtained taking the
first line and the first column and letting the other row and
column vary, i.e.
$$
{\tilde A_{i\,j}}:=\det A_{1\,i}^{1\,j}\,\,\,,\quad 2\leq i,j\leq
g.
$$
We observe that decomposing the matrix $A$ in blocks $$A=\pmatrix
a_{11}&^t z\\ w&B\endpmatrix,$$ with $B$ a $(g-1)\times( g-1)$
matrix, and $z, w\in \C^{g-1}$, we have
$$
{\tilde A_{i\,j}}:=\det (a_{11}B-w^t z)
$$
With these notations we have

\begin{lm}
\begin{equation}
\label{det1} a_{1\,1}^{g-2}\det A= \det {\tilde A}
\end{equation}
\end{lm}
\begin{proof} It is trivial when $a_{1\,1}=0$ and it is
 an immediate consequence of $$A=\pmatrix 1&0\\  w/a_{11}&
{1\over a_{11}}I\endpmatrix \pmatrix a_{11}&^t z\\ 0&{\tilde
A}\endpmatrix$$ when $a_{1\,1}\neq 0$
\end{proof}

We denote $N:= g(g+1)/2$ and for any $v\in\C^g$ we denote by $v^2$
its symmetric tensor square. Then the following is true
\begin{lm}
Let $v_1, \dots,v_N\in\C^g$. Then
$$
(N!) v_{1} ^2\wedge v_{2} ^2\wedge \dots\wedge v_{N}^2=\sum_{s\in
S_N}{\rm sign}(s) (v_{s(1)}\wedge v_{s(2)}\wedge\dots\wedge
v_{s(g)})\cdot
$$
$$
(v_{s(1)}\wedge v_{s(g+1)}\wedge v_{s(g+2)}\wedge\dots\wedge
v_{s(2g-1)})\cdot (v_{s(2)}\wedge v_{s(g+1)}\wedge v_{s(2g
)}\wedge\dots\wedge v_{s(3g-3)})
$$
$$
(v_{s(3)}\wedge v_{s(g+2)}\wedge v_{s(2g)}\wedge\dots\wedge
v_{s(4g-6)})\dots (v_{s(g)}\wedge v_{s(2g-1)}\wedge v_{s(3g-3
)}\wedge\dots\wedge v_{s(N)})
$$
\end{lm}
\begin{proof}
Since the LHS is $SL(g,\C)$ invariant, it can be expressed as a
polynomial in determinants of $g\times g$ minors of the $g\times
N$ matrix with columns being $v$'s. Moreover, this polynomial must
be homogeneous of degree $g+1$ in these determinants, each $v_i$
has to appear in it exactly twice, and it has to be skew-symmetric
in $v$'s. For this reason first we sum over all possible
permutations with the signs. We further observe that if the same
two vectors appear in two different determinants, then the
expression vanishes. Thus the expression has to be a sum of
monomials each of degree $g+1$ in the determinants, such that each
$v$ appears exactly twice, and no pair of $v$'s appears twice in
two different determinants. Thus the expression is forced to be
exactly that of the statement, up to a multiplicative constant,
which is easily computed.
\end{proof}

\section{$\t$'s and $\T$'s}
A special case of Riemann's bilinear addition theorem for theta
functions (see \cite{ig1},\cite{farkasrauch},\cite{mumford}) is
\begin{equation}
\label{tT} \T[\a](\tau,z)\T[\a+\e](\tau,0)
=\frac{1}{2^g}\sum\limits_{\s\in(\Z/2\Z)^g}(-1)^{ \a\cdot\s}\tt\e{
\s}( \tau, z)\tt\e\s (\tau,z)
\end{equation}
which is valid for all $\tau$ and $z$. Taking a sum of these with
different signs, we get, for any $\de\in(\Z/2\Z)^g$
\begin{equation}
\label{tT2}
\begin{matrix}
\sum\limits_{\a\in(\Z/2\Z)^g}
(-1)^{\a\cdot\de}\T[\a](\tau,z)\T[\a+\e](\tau,0) =\\
\frac{1}{2^g}\sum\limits_{\a,\,\s\in(\Z/2\Z)^g}
(-1)^{\a\cdot(\s+\de)}\tt\e{\s}(\tau,z)\tt\e\s(\tau,z)=\tt\e{\de}
(\tau, z)\tt\e\de(\tau,z).
\end{matrix}
\end{equation}

We assume that the characteristic $[\e,\de]$ is odd, differentiate
this relation with respect to $z_i$ and $z_j$, and then evaluate
at $z=0$. Denoting by ${\bf C}_{\e\,\de}(\tau)$ the $g\times g$
symmetric matrix with entries
$$
{\bf C}_{\e\,\de,\,\, ij}( \tau) :=2\p_{z_i}\tt\e{\de}( \tau,
0)\p_{z_j}\tt\e\de (\tau,0),
$$
and by $\bf A_{\e\,\de}(\tau)$ --- the $g\times g$ symmetric
matrix with entries
$$
{\bf A}_{\e\,\de,\,\, ij}(\tau):=
(\p_{z_i}\p_{z_j}\T[\de](\tau,0))\T[\e](\tau,0)-
(\p_{z_i}\p_{z_j}\T[\e](\tau,0))\T[\de](\tau,0),
$$
we thus have (see \cite{gsm}) --- notice that ${\bf
C}_{\e\,\de}=0$ if $[\e,\,\de]$ is even
\begin{lm}
\label{AC}
\begin{equation}
{\bf C}_{ \e\,\de}( \tau) =\frac{1}{2}
\sum\limits_{\a\in(\Z/2\Z)^g}(-1)^{\a\cdot \de }{\bf A}_{
{\e+\a}\, \a}( \tau)
\end{equation}
\end{lm}
and the ``inverse''
\begin{lm}
\begin{equation}
\label{Adef} {\bf A}_{ {\e+\a}\,\a}( \tau) =
\frac{1}{2^{g-1}}\sum\limits_{\b\in(\Z/2\Z)^g}(-1)^{ \a\cdot \b
}{\bf C}_{\e\,\b}( \tau).
\end{equation}
\end{lm}

We remark also that
\begin{equation}
{\bf C}_{ \e\,\de}( \tau)=2grad\,\,\tt\e\de(\tau)\, \,^{t}
grad\tt\e\de(\tau).
\end{equation}

\section{Generalized Jacobi's derivative formulas}
To generalize the first result of the introduction, we introduce
the matrix-valued differential operator
$$
\cal D:= \left(\begin{array}{rrrr}
\,\frac{\p}{\p\tau_{11}}&\frac{1}{2}\frac{\p}
{\p\tau_{12}}&\dots&\frac{1}{2}\frac{\p}{\p\tau_{1 g}}\\
\frac{1}{2}\frac{\p}{\p \tau_{21}}&\frac{\p}{\p
\tau_{22}}&\dots&\frac{1}{2}\frac{\p}{\p\tau_{2 g}}\\
\dots&\dots&\dots&\dots\\
\frac{1}{2}\frac{\p}{\p \tau_{g 1}}& \dots&\dots& \,\
\frac{\p}{\p\tau_{g g}}\end{array}\right).
$$
Then we have
\begin{thm}[First generalization]
For any $\e\ne\de$ the following holds:
\begin{equation}
\label{Eq}
\begin{matrix}
{\rm c}\T[\de]^{2g} \det(\cal D(\T[\e]/\T[\de])\quad\hfill\\
\hfill=\!\!\!\!\!\!\!\!\sum\limits_{\lbrace\a_{i_1},\dots,
\a_{i_g}| [\e+\de,\a_{i_j}] \,{\rm
odd}\rbrace}\!\!\!\!\!\!\!\!(-1)^{ \de\cdot
(\a_{i_1}+\dots+\a_{i_g})} D([\e+\de,\,\a_{i_1}],\dots
[\e+\de,\,\a_{i_g}])^2
\end{matrix}
\end{equation}
for some computable constant ${\rm c}$.
\end{thm}
\begin{proof}
For any characteristics $\e,\de$ we have by definition
$$
\T[\de]^2\cal D(\T[\e]/\T[\de])(\tau)=4\pi i{\bf A}_{ \e
\de}(\tau).
$$
Thus
$$
\T[\de]^{2g} \det\,(\cal D(\T[\e]/\T[\de]))(\tau)=(4\pi i)^g
\det({\bf A}_{ \e \de})(\tau).
$$
Now, using the result of lemma 2, we get
$$
\det{\bf A}_{\e\de}(\tau)=\det\left(\frac{1}{2^{g-1}}\sum_{\lbrace
\a|[\e+\de,\,\,\a]\,{\rm odd}\rbrace}(-1)^{\de\cdot \a}{\bf C}_{
\e+\de,\,\,\a})(\tau)\right)=
$$
$$
\det\left(\frac{1}{2^{g-1}}\sum_{\lbrace\a|[\e+\de,\,\,\a]\,{\rm
odd}\rbrace}(-1)^{\de\cdot\a}\,\,grad\,\tt{\e+\de}\a(\tau)^{t}\,\,
grad\,\tt{\e+\de}\a(\tau)\right).
$$
When we expand this determinant, each summand will be of the type
$$
{\rm sign}(\mu)\left((-1)^{\de\cdot \a_1 }\p_{\mu(1)}
\tt{\e+\de}{\a_1}( \tau)\cdots (-1)^{\de\cdot
{\a_g}}\p_{\mu(g)}\tt{\e+\de}{\a_g} ( \tau)\right)
$$
$$
\quad{\rm sign}(\s)\left(\p_{\s(1)}\tt{\e+\de}
{\a_1}(\tau)\cdots\p_{\s(g)}\tt{\e+\de}{ \a_n} (\tau)\right)
$$
for some permutations $\s$ and $\mu$. Taking the sum of these for
all possible permutations $\s$ and $\mu$ gives exactly the square
of the jacobian determinant, so that we end up with
$$
\left(\frac{\pi^2}{2^{g-2}}\right)^g
\sum\limits_{\a_{i_1},\dots,\a_{i_g}\in
(\Z/2)^g}(-1)^{\de\cdot(\a_{i_1}+\dots+\a_{i_g})}D([\e+\de,\,\a_{
i_1}] ,\dots[\e+\de,\,\a_{i_g}])^2,
$$
proving the theorem.
\end{proof}

At this point, we observe that our relations are not trivial. In
fact each term appearing in the RHS is not identically zero,
\cite{sm3}. We remark that the set of characteristics appearing in
the jacobian determinant above is syzygetic, while in all the
other generalizations of Jacobi's derivative formula only azygetic
sets appear, cf. \cite{ig2}.

If we would like to have relations involving the derivatives of
the second order theta constants with respect to $\tau_{ij}$, then
since the matrix ${\bf C}_{\e\,\de}(\tau)$ has rank $1$, we have
the following
\begin{prop}
If $g\geq 2$, then
\begin{equation}
\label{cor}
\det\left(\sum\limits_{\a\in(\Z/2\Z)^g}(-1)^{\a\cdot\de}{\bf
A}_{{\e+\a}\, \a}(\tau)\right)=0.
\end{equation}
\end{prop}

For the LHS of (\ref{Eq}), we have the following non-vanishing:
\begin{thm}
For all possible pairs $\e\not=\de\in(\Z/2\Z)^g$ the expression
$$
\T[\de]^{2g} \det\,(\cal D(\T[\e]/\T[\de]))
$$
is not identically zero in $\tau$.
\end{thm}
\begin{proof}
We shall prove a slightly more general result: that for any pair
of distinct even characteristics $[\e,\a]$ and $[\de,\b]$ the
expression
$$
\tt\de\b^{2g} \det\,\left(\cal D\left(\tt\e\a
/\tt\de\b\right)\right)(\tau)
$$
is not identically zero.

Indeed, we know that the symplectic group $\Gamma _g$ acts doubly
transitively on the set of even characteristics, and we have the
following transformation formula
$$
\begin{matrix}
\cal D\left(\theta \bmatrix M\pmatrix \e\cr
\a\endpmatrix\endbmatrix\ /\ \theta\bmatrix M\pmatrix\de\cr
\b\endpmatrix\endbmatrix\right)(M\cdot\tau )\\
= \phi(\e,\a,\de, \b,M)(C\tau+D)^t\cal D\left(\tt\e\a
/\tt\de\b\right)(\tau)(C\tau+D ).
\end{matrix}
$$
with $\phi(\e,\a,\de,\b,M)$ an eighth root of unity.

We learnt from \cite{bz} that for some specific $[\e,\,0]$ and
$[\de,\,0]$ we have $\det\cal D\left(\tt\e0 /\tt\de0\right)$ not
identically zero; since all such expressions are permuted by the
symplectic group action, they are all not identically zero.
\end{proof}

\medskip
For the second generalization of Jacobi's derivative formula, for
any set of $N+1:=\frac{1}{2}g(g+1)+1$ characteristics $\e_0, \e_1,\dots,\e_N$ we
introduce the matrix
$$
M(\e_0, \e_1,\dots,\e_N)(\tau):=\pmatrix
\hskip6.5mm\T[\e_0]&\dots&\hskip6.5mm\T[\e_N]\\
\frac{\p}{\p\tau_{11}}\T[\e_0]&\dots&\frac{\p}{\p\tau_{11}}\T[\e_N]\\
\frac{\p}{\p\tau_{12}}\T[\e_0]&\dots&\frac{\p}{\p\tau_{12}}\T[\e_N]\\
\dots&\dots&\dots\\
\frac{\p}{\p\tau_{gg}}\T[\e_0]&\dots&\frac{\p}{\p\tau_{gg}}\T[\e_N]
\endpmatrix
$$
It is a well-known fact (see \cite{sasaki}) that $\det M(\e_0,
\e_1,\dots,\e_N)(\tau)$ is a modular form of weight
$\frac{1}{4}(g+2)(g+3)$ relatively to $\Gamma_g (2,4)$. Moreover,
if it is not identically zero, then the theta constants
$$
\T[\e_0], \T[\e_1],\dots ,\T[\e_N]
$$
are algebraically independent. We set
$$
\de_k:=\e_0+\e_k
$$
and let $\cal C_{\e\,\b }(\tau)$ be the vector in $\C^N$ with
entries $C_{\e\,\b,\,\, ij}(\tau)$ (before we thought of $C$ as a
matrix, but now we write down all the matrix elements in a single
vector). Using the results of Lemmata 4, 2 and 1, we get the
following
\begin{prop}[Second generalization]
For some computable constant $c$
$$ c\T[\e_0]^{N-1}\det M(\e_0,
\e_1,\dots,\e_N)(\tau)
$$
$$
=\sum \limits_{\b_1,\ldots, \b_N\in (\Z/2\Z)^g} (-1)^{\e_k \b_k}
\det( \cal C_{\de_1\,\b_1}\wedge \cal
C_{\de_2\,\b_2}\wedge\dots\wedge \cal C_{\de_N\,\b_N})
$$
By lemma 2 the RHS can be expressed as a homogeneous polynomial of
degree $g+1$ in jacobian determinants.
\end{prop}
\begin{rem}
Recalling the definition of theta functions, we can expand all of
the above identities in Fourier series in $\tau$ or equivalently
in power series in $q_{ij}:=\exp \tau_{ij}$. The coefficients of
these expansions will then be some rather complicated but quite
natural combinatorial quantities in several variables, and the
equality of the RHS and LHS of any of the above would then yield a
non-trivial multidimensional combinatorial identity, which it
would be interesting to understand and prove combinatorially.
\end{rem}

\section{An application in genus 2}
We will now work out in detail the situation in the case of genus
2. Indeed let us write down (\ref{Eq}) for $\e=[00]$ and
$\de=[10]$:
\begin{equation}
\label{gen2} \T[10]^ 4 \det(\cal D
(\T[00]/\T[10])=cD([10,10],[10,11])^2,
\end{equation}
with $c$ a known constant.

Using lemma \ref{AC} to express the RHS in terms of theta
constants of the second order and their derivatives, we get (we
denote $\p_{ij}:=\p_{\tau_{ij}}$ to simplify notations)
$$
\quad\left(\T[00]\p_{11}\T[10]-\T[10]\p_{11}\T[00]\right)
\left(\T[00]\p_{22}\T[10]-\T[10]\p_{22}\T[00]\right)
$$
$$
-\left(\T[00]\p_{12}\T[10]-\T[10]\p_{12}\T[00]\right)^2\quad\quad
$$
$$
+\left(\T[01]\p_{11}\T[11]-\T[11]\p_{11}\T[01]\right)
\left(\T[01]\p_{22}\T[11]-\T[11]\p_{22}\T[01]\right)
$$
$$
-\left(\T[01]\p_{12}\T[11]-\T[11]\p_{12}\T[01]\right)^2=0
$$
Clearly we get the same equation if we chose $\e=[01]$ and
$\de=[11]$. Thus in all we get three different equations.

In \cite{bz} it is shown that in genus 2 there are $2\cdot
2^2+2=10$ algebraically independent quantities among 4 theta
constants of the second order and their $4\cdot 3=12$ derivatives.
Thus there are 6 non-trivial algebraic relations among theta
constants and their first-order derivatives. So far we have
obtained three such equations, and three more can be obtained by
writing down formula (\ref{cor}).

In \cite{bz} some other 6 independent relations are given. We
shall prove that they are all consequences of (\ref{Eq}) and
(\ref{jac1}).

Indeed, for genus two formula (\ref{jac1}) reads
$$
D([10,10],[10,11])^2=\left(\tt {11} {00} \tt {11} {11}\tt {01} {00}
\tt {01} {10}\right)^2 .
$$
Applying (\ref{tT2}) to rewrite the RHS in terms of theta
constants of the second order, we finally see that
$$
\T[10]^ 4 \det(\cal D (\T[00]/\T[10])
$$
is a polynomial in the theta constants of the second order. This
equation is up to a rational function  equal to one of the
equations in \cite{bz}. From the other choices of characteristics
$\e$ and $\de$ we get the other 5 equations. We observe that these
5 equations can also be obtained from the first one by the action
of $\Gamma_2$. In this spirit we remark that in genus 1, we have
four variables and the relation is an immediate consequence of
(\ref{jac0}) and (\ref{jac4}). These relations can also be
obtained considering the determinant of (\ref{cor}).

This method allows us to give a conjectural description of the
situation in the genus 3 case. By the results of \cite{bz} we know
that among 56 variables (all $\T[\e] $ and their derivatives
$\p_{ij}\T[\e]$) there are 21 algebraically independent ones. Thus
there are 35 algebraic relations. We know that there is a unique
polynomial relation among the $\T[\e]$, of degree 16, cf.
\cite{vg}; let us denote it by $R(\tau)$. Thus we have
$$
R(\tau)=0, \quad \p_{ij}R(\tau)=0
$$
and other 28 relations obtained by applying (\ref{cor}), since in
genus 3 there are exactly 28 odd characteristics. So in total we
have 35 relations that we conjecture to be algebraically
independent.

\section{Characterization of the reducible locus}
We finish the paper by giving a characterization of the locus of
reducible abelian varieties. Different characterizations of the
reducible locus are known: in \cite{el} it is characterized in terms
of the dimension of the singular locus of the theta divisor, in
\cite{sasaki} --- in terms of the non-maximality of the rank of
matrix $P(\tau)$ with rows $( \T[\e],\frac{\p}
{\p\tau_{ij}}\T[\e])$, and columns corresponding to all
$\e\in(\Z/2\Z)^g$, in \cite{sm1} --- in terms of the vanishing of
certain theta constants. Here we use the vanishing of certain first
derivatives of theta functions evaluated at zero.
\begin{prop}
A ppav with a period matrix $\tau$ is reducible if and only if
there exist some $M\in\Gamma_g$ and some $k<g$ such that if we
write any odd characteristic $[\e,\de]$ as
$[\e_1\,\e_2,\de_1\,\de_2]$, where $[\e_1,\de_1]$ is a
$k$-dimensional characteristic, and $[\e_2,\de_2]$ is
$(g-k)$-dimensional, then
$$
\p_{z_i}\tt\e\de(M\cdot\tau,z)|_{z=0}=0
$$
for all $i\leq k$ for $[\e_1 ,\de_1]$ even, and for all $i>k$ for
$[\e_1,\de_1]$ odd.
\end{prop}
\begin{proof}
Suppose the period matrix $M\cdot\tau$ splits as $M\cdot\tau=
\pmatrix \tau_1&0\cr 0&\tau_2\endpmatrix $, with $\tau_1\in\H_k$
and $\tau_2\in\H_{g-k}$, so that the theta functions with
characteristics factor as follows:
$$
\tt\e\de(M\cdot\tau,z)=\tt{\e_1}{\de_1}(\tau_1,z_1)\cdot
\tt{\e_2}{\de_2}(\tau_2,z_2).
$$
The vanishing of the derivatives in question is immediate for
$M\cdot\tau$ by differentiating and evaluating at $z=0$; thus the
``only if'' part is proven. For the ``if'' part, assume the
vanishing of derivatives as stated. Then according to (\ref{Adef})
we have ${\bf C}_{\e\de,ij}=0$ and consequently ${\bf
A}_{\e\de,ij}=0$ for $1\leq i\leq k<j\leq g$ and all odd
$[\e,\de]$. Thus the matrix $P(\tau)$ does not have maximal rank,
and thus, by the results of \cite{sasaki}, it corresponds to a
reducible abelian variety.
\end{proof}

\end{document}